\documentclass[amssymb, 11pt]{amsart}
\usepackage{latexsym}

\newcommand{\ee}{\mathbb{E}}
\newcommand{\pp}{\mathbb{P}}
\newcommand{\ii}{\mathbb{I}}
\newcommand{\ex}{\mathrm{Exp}(1)}

\newcommand{\rr}{\mathbb{R}}
\newcommand{\fp}{f^\prime}
\newcommand{\fpp}{f^{\prime\prime}}
\newcommand{\ww}{W^\prime}
\newcommand{\frc}[2]{{\textstyle{\frac{#1}{#2}}}}
\def\norm#1{\left\Vert#1\right\Vert}
\def\abs#1{\left\vert#1\right\vert}
\def\law{\mathcal{L}}
\def\IE{\mathbb{E}}
\def\IP{\mathbb{P}}
\def\IR{\mathbb{R}}
\def\II{\mathbb{I}}
\def\Exp{\mathrm{Exp}}

\def\ahalf{{\textstyle\frac{1}{2}}}

\def\tsfrac#1#2{{\textstyle\frac{#1}{#2}}}

\newtheorem{lmm}{Lemma}

\newlength{\standardunitlength}
\newtheorem{prop}{Proposition}[section]

\newtheorem{lemma}[prop]{Lemma}

\newtheorem{theorem}[prop]{Theorem}

\begin{document}

\begin{center}
{\bf Exponential Approximation by Stein's Method and Spectral Graph
Theory}
\end{center}

\begin{center}
Running head: Exponential Approximation by Stein's Method
\end{center}

\begin{center}
Version of 8/16/08
\end{center}

\begin{center}
By Sourav Chatterjee, Jason Fulman, and Adrian R\"{o}llin
\end{center}

\author{Sourav Chatterjee}
\address{Department of Statistics\\
University of California, Berkeley\\
Berkeley, CA 94720-3860}
\email{sourav@stat.berkeley.edu}

\author{Jason Fulman}
\address{Department of Mathematics\\
University of Southern California\\
Los Angeles, CA 90089} \email{fulman@usc.edu}

\author{Adrian R\"{o}llin}
\address{Department of Mathematics\\
National University of Singapore\\
Singapore 117543} \email{matar@nus.edu.sg}

{\bf Abstract}: General Berry-Ess\'{e}en bounds are developed for the
exponential distribution using Stein's method. As an application, a sharp
error term is obtained for Hora's result that the spectrum of the
Bernoulli-Laplace Markov chain has an exponential limit. This is the first
use of Stein's method to study the spectrum of a graph with a non-normal
limit.

\begin{center}
\end{center}

\begin{center}
\end{center}

\section{Introduction}

    This paper develops general Berry-Ess\'{e}en bounds for the exponential distribution
using Stein's method. Two of our main results are given by the following
statements. We let $\ii[A]$ denote the indicator function of an event $A$.

\begin{theorem} \label{atheorem1} Assume
that $W$ and $W'$ are non-negative random variables on the same
probability space such that $\law(W')=\law(W)$. Then, if $Z\sim\Exp(1)$,
we have for any $t>0$ and any constant $\lambda>0$
\begin{eqnarray*}
    \abs{\IP[W\leq t]-\IP[Z\leq t]}
 & \leq & \IE\abs{(\lambda^{-1}\IE(D|W)+1)\ii[W>0]} +
\IE\abs{{\textstyle\frac{1}{2\lambda}} \IE(D^2|W)-1}\\
    & & +{\textstyle\frac{1}{6\lambda}}\IE\abs{D}^3
    + {\textstyle\frac{1}{2\lambda}}
        \IE\left(D^2\ii[\abs{W-t}\leq\abs{D}]\right),
\end{eqnarray*}
where $D := W' - W$.
\end{theorem}

\begin{theorem} \label{basicexp} Assume
that $W$ and $W'$ are non-negative random variables on the same
probability space such that $\law(W')=\law(W)$ and
\begin{equation*}
    \IE(D|W) = -\lambda(W-1),
\end{equation*}
where $\lambda>0$ is a fixed constant. Then if $Z\sim\Exp(1)$, we have for
any $t>0$,
\begin{equation*}
\begin{split}
\left|\IP[W\leq t] - \IP[Z\leq t]\right|  &  \leq\frac{\IE\left|2\lambda W
- \IE(D^2|W) \right|}{2\lambda t} +
\frac{\IE|D|^3\max\{t^{-1},2t^{-2}\}}{6\lambda}\\ &\quad +
\frac{\IE\left\{D^2 \ii \left[|W-t|\leq|D|\right]\right\}}{\lambda t}.
\end{split}
\end{equation*}
where $D:=W'-W$.
\end{theorem} The use of a pair $(W,W')$ is similar to the exchangeable
pairs approach of Stein for normal approximation \cite{St1}, but in the
spirit of \cite{Ro}, throughout this paper we require only the weaker
assumption that $W$ and $W'$ have the same law. It can be challenging to
obtain good bounds on the error terms in Theorems \ref{atheorem1} and
\ref{basicexp}, and we also develop a number of tools for doing that.

Before continuing, we mention that this is not the first paper to study
exponential approximation by Stein's method. Indeed, earlier works, in the
more general context of chi-squared approximation, include Mann
\cite{Man}, Luk \cite{Lu}, and Reinert \cite{Re} (which also includes a
discussion of unpublished work of Pickett). The paper \cite{Man} uses
exchangeable pairs, whereas \cite{Lu} and \cite{Re} use the generator
approach to Stein's method. However all of these papers focus on
approximating expectations of smooth functions of $W$, rather than
indicator functions of intervals, and so do not give Berry-Ess\'{e}en
theorems. Moreover, the examples in \cite{Lu} and \cite{Re} are about sums
of independent random variables, whereas our example involves dependence.

    Our main example is the spectrum of the
    Bernoulli-Laplace Markov chain. This Markov chain was
    suggested as a model of diffusion and has the following
    description. Let $n$ be even. There are two urns, the first
    containing $\frac{n}{2}$ white balls, and the second
    containing $\frac{n}{2}$ black balls. At each stage, a ball is
    picked at random from each urn and the two are
    switched. Diaconis and Shahshahani \cite{DS} proved that
    $\frac{n}{8}\log(n) + \frac{cn}{2}$ steps suffice for this
    process to reach equilibrium, in the sense that the total
    variation distance to the stationary distribution is at most
    $ae^{-dc}$ for positive universal constants $a$ and $d$. In
    order to prove this, they used the fact that the spectrum of
    the Markov chain consists of the numbers
    $1-\frac{i(n-i+1)}{(n/2)^2}$ occurring with multiplicity ${n
    \choose i} - {n \choose i-1}$ for $1 \leq i \leq \frac{n}{2}$
    and multiplicity 1 if $i=0$. Hora proved the following result,
    which shows that the spectrum of the Bernoulli-Laplace chain
    has an exponential limit.

\begin{theorem} \label{hoexp} (\cite{Ho1}) Consider the uniform measure
on the set of the ${n \choose \frac{n}{2}}$
 eigenvalues of the Bernoulli-Laplace Markov chain. Let $\tau$ be a
 random eigenvalue chosen from this measure. Then as $n \rightarrow
 \infty$, the random variable $W:= \frac{n}{2} \tau+1$ converges
 in distribution to an exponential random variable with mean $1$. \end{theorem}

    As an application of our general Berry-Ess\'{e}en bound, the
    following result will be proved.

\begin{theorem}
 \label{sharper} Let $Z \sim \ex$, and let $W$ be as in Theorem
 \ref{hoexp}. Then \[ |\pp \{W \leq t \} -\pp \{Z \leq
 t \}| \leq \frac{C}{\sqrt{n}}, \] for all $t$, where $C$ is a universal constant.
Moreover this rate is sharp in the sense that there is a sequence of $n$'s
tending to infinity, and corresponding $t_n$'s such that \[ |\pp(W \leq
t_n) - \pp(Z \leq t_n)| = \frac{2e^{-2}}{\sqrt{n}} + O(1/n).\]
\end{theorem}

Note that the Bernoulli-Laplace Markov chain is equivalent to random walk
on the Johnson graph $J(n,k)$ where $k=\frac{n}{2}$. The vertices consist
of all size $k$ subsets of $\{1,\cdots,n\}$, and two subsets are connected
by an edge if they differ in exactly one element. From a given vertex,
random walk on the Johnson graph picks a neighbor uniformly at random, and
moves there.

One reason why our method for proving Theorem \ref{sharper} is of interest
(despite the existence of a more elementary argument for a weaker version
of Theorem \ref{sharper} sketched at the end of Section \ref{n/2}) is that
Theorem \ref{sharper} is in fact a small piece of a much larger program.
To explain, limit theorems for graph spectra (especially Cayley graphs and
finite symmetric spaces) have been studied by many authors and from
various perspectives; some references are \cite{Ho1}, \cite{Ho2},
\cite{Ke}, \cite{F1}, \cite{F2}, \cite{F3}, \cite{F4}, \cite{Sn},
\cite{ShSu}, \cite{T1}, and \cite{T2}. In particular, the references
\cite{T1} and \cite{T2} describe some challenging conjectures where the
limit distribution is the semicircle law and relate them to deep work in
number theory. With the long-term goal of making progress on these
conjectures, the paper \cite{F4} gave some general constructions for
applying Stein's method to study graph spectra, and worked out examples in
the case of normal approximation. The current paper works out an
exponential example, and is excellent evidence that these constructions
will prove useful in other settings where the spectrum has a non-normal
limit. We also emphasize that while there are papers such as \cite{GT}
which obtain non-normal limit theorems with an error term in spectral
problems, they study the spectrum of random objects, whereas our work, and
the conjectures of \cite{T1}, \cite{T2}, all pertain to the spectrum of a
sequence of fixed, non-random graphs.

We also mention that an additional reason for studying the spectrum of the
Bernoulli-Laplace chain is that it is closely related to the spectrum of
the random transposition walk (and so with representation theory of the
symmetric group). Indeed, from \cite{Sc}, the eigenvalues of the
Bernoulli-Laplace chain can be expressed as $\frac{2(n-1)}{n} \mu -
\frac{(n-2)}{n}$ as $\mu$ ranges over a subset of eigenvalues of the
random transposition walk. This relationship is not surprising given that
the Bernoulli-Laplace chain transposes balls from different urns at each
step. But together with the large body of work on Kerov's central limit
theorem for the spectrum of the random transposition walk (\cite{Ke},
\cite{F2}, \cite{F3}, \cite{F4}, \cite{Sn}, \cite{IO}, \cite{Ho2}), it
does make the problems studied in the current paper very natural. As a
final justification for the current paper, we believe that the example in
it will serve as a useful testing ground for other researchers in Stein's
method (certainly it helped us in developing our Berry-Ess\'{e}en
theorems).

The organization of this paper is as follows. Section \ref{v1} proves our
first general Berry-Ess\'{e}en bound for the exponential law, namely
Theorem \ref{atheorem1} above, and develops tools for analyzing the error
terms which appear in it. Section \ref{v2} proves our second general
Berry-Ess\'{e}en bound for the exponential law, namely Theorem
\ref{basicexp} above, and develops tools paralleling those in Section
\ref{v1} for analyzing the error terms. Section \ref{n/2} treats our main
example (spectrum of the Bernoulli-Laplace chain), proving Theorem
\ref{sharper}. An interesting feature of the proof is that it uses theory
from both of Sections \ref{v1} and \ref{v2}, to treat the cases of small
and large $t$ respectively. Finally, Appendix \ref{Alg} gives an algebraic
approach to the exchangeable pair and moment computations in Section
\ref{n/2}, linking it with the constructions of \cite{F4}. This is not
essential to the proofs of any of the results in the main body of the
paper, but does motivate the exchangeable pairs used in the paper, which
could be difficult to guess.

\section{Berry-Ess\'{e}en Bound for the Exponential Law: Version 1} \label{v1}

A main purpose of this section is to prove Theorem \ref{atheorem1} from
the introduction, and to develop tools for analyzing the error terms which
arise in it. To begin we make some remarks concerning the statement of
Theorem \ref{atheorem1}.

{\it Remarks:}
\begin{enumerate}

\item In our main example (see Section \ref{n/2}), the relation
$\ee(D|W)=-\lambda$ is satisfied for all $W>0$. Hence the first error term
in Theorem \ref{atheorem1} will vanish. In the spirit of \cite{RR}, one
could also have that $\ee(D|W)=-\lambda + R$ for some non-trivial random
variable $R$.

\item Although $W$ is allowed to attain the value $0$ (and does, in our main
example), the conditional expectation $\ee(D|W=0)$ (i.e. the ``drift'' at
zero) does not enter in the first term of the bound.

\item In our main example (see Section \ref{n/2}), $\ee(D^2|W)=2 \lambda$
and so the second error term also vanishes. The third error term
$\ee|D|^3$ can be bounded using the Cauchy-Schwarz inequality $\ee|D|^3
\leq \sqrt{\ee|D|^2 \ee|D|^4}$. The error term which is difficult to bound
in practice is the fourth error term, and later in this section we develop
suitable tools (see Theorems \ref{tool1} and \ref{tool2}).

\end{enumerate}

Before embarking on the proof of Theorem \ref{atheorem1}, we recall the
main idea of Stein's method in our context. As observed by Stein
\cite{St2}, a random variable $Z$ on $[0,\infty)$ is $\Exp(1)$ if and only
if $\ee[f'(Z)-f(Z)]= -f(0^+)$ for all functions $f$ in a large class of
functions (whose precise definition we do not need). Here $f(0^+)$ is the
limiting value of $f(a)$ as $a$ approaches $0$ from the right. Stein's
characterization of the exponential distribution motivates the study of
the function $f(x)$ solving the equation
\[ f'(x) - f(x) = I[x \leq t] - (1-e^{-t}), \qquad x\geq 0.\] Indeed, for
this $f$ one has that \[ \pp(W \leq t) - \pp(Z \leq t) = \ee[f'(W)-f(W)],
\] and the problem becomes that of bounding $\ee[f'(W)-f(W)]$.

We begin with the following lemma.

\begin{lemma} \label{boundsol} For every $t> 0$, the function
\begin{equation} \label{first}
    f(x) := e^{-(t-x)^+} - e^{-t},\qquad x\geq 0, \end{equation}
(where in \eqref{first} we define the derivative $f'(t) := f'(t^-)$),
satisfies the differential equation
\begin{equation}\label{steinequation}
    f'(x) - f(x) = \ii[x \leq t] - (1-e^{-t}), \qquad x\geq 0,
\end{equation}
and the bounds
\begin{equation} \label{boundsonsolution}
    \norm{f}_\infty \leq 1,\quad \norm{f'}_\infty \leq 1, \quad
    \sup_{x,y\geq0}\abs{f'(x)-f'(y)}\leq 1.
\end{equation}
The second derivative  $f''$, defined for every $x\neq t$, satisfies
\begin{equation} \label{boundsecondderivative}
    \sup_{x\neq t} \abs{f''(x)} \leq 1.
\end{equation}
\end{lemma}

\begin{proof} Write
\[
    f(x) = \begin{cases}
           e^{-t+x}-e^{-t} & \text{if $x\leq t$,}\\
           1-e^{-t}         & \text{if $x> t$.}
           \end{cases}
\]
Together with the definition of $f'(t)$ this yields
\[
    f'(x) = \begin{cases}
           e^{-t+x} & \text{if $x\leq t$,}\\
           0         & \text{if $x> t$.}
           \end{cases}
\]
Thus, on $x\leq t$, \[    f'(x) - f(x)  = 1 - (1-e^{-t})
\]
which is \eqref{steinequation}, and on $x> t$
\[
    f'(x) - f(x) = 0 - (1-e^{-t})
\]
which again is \eqref{steinequation}. The bounds \eqref{boundsonsolution}
and \eqref{boundsecondderivative} are straightforward; to obtain the last
bound in \eqref{boundsonsolution} note that $f'$ is non-negative.
\end{proof}

Now we give a proof of Theorem \ref{atheorem1}.

{\sc Proof of Theorem \ref{atheorem1}} Using \eqref{steinequation} it is
clear that we only need to bound $\ee(f'(W)-f(W))$. Fix $t>0$ and let
$F(x) := \int_0^x f(y)dy$. By Taylor expansion,
\begin{equation}\begin{split}\label{amainequation}
 0 & = \ee\left(F(W')-F(W)\right) \\
   & = \ee\left(D f(W)\right) +
        \ee\left(D^2\int_0^1(1-s)f'(W+sD)ds\right)\\
   & = \ee\left(D f(W) + \ahalf D^2 f'(W)\right) + \ee\left(D^2 J\right)\\
\end{split}\end{equation}
where
\[
    J := \int_0^1 (1-s)(f'(W+sD)-f'(W))ds.
\]

Let $A$ be the event that $W\wedge W' \leq t \leq W \vee W'$. On $A^c$ we
thus have for every $0\leq s \leq 1$
\begin{equation}    \label{aboundsproof1}
    \abs{f'(W+sD)-f'(W)} \leq s\abs{D}
\end{equation}
by \eqref{boundsecondderivative},  whereas on $A$ we have
\begin{equation}    \label{aboundsproof2}
    \abs{f'(W+sD)-f'(W)} \leq 1
\end{equation}
by \eqref{boundsonsolution}. Dividing \eqref{amainequation} by $\lambda$
and noting that $f(0) = 0$, and thus $f(W)=\ii[W>0]f(W)$, we can use this
to obtain that
\begin{equation} \begin{split} \label{amainequation2}
    \ee(f'(W)-f(W)) &
    =  \ee(f'(W)-f(W)) - \frac{1}{\lambda} \left(\ee\left(F(W')-F(W)\right)\right)\\
    &= -\ee\left((\tsfrac{1}{\lambda}\ee(D|W)+1)f(W)\ii[W>0]\right) \\
    &  + \ee\left((1-\tsfrac{1}{2\lambda}\ee(D^2|W)) f'(W)\right) \\
    &  - \tsfrac{1}{\lambda}\ee\left(\ii[A^c]D^2 J\right)
    - \tsfrac{1}{\lambda}\ee\left(\ii[A]D^2 J\right).\\ \end{split}
\end{equation} Invoking the bounds \eqref{aboundsproof1} and
\eqref{aboundsproof2}, we have
\[
    \ii[A^c]D^2\abs{J} \leq {\textstyle\frac{1}{6}}|D|^3,\\
    \ii[A]D^2\abs{J} \leq \ahalf D^2 \ii[\abs{W-t}\leq |D|],
\]
where the second inequality uses the fact that  $A$ implies $\abs{W-t}
\leq |D|$. Combining these bounds with \eqref{amainequation2} and the
bounds $||f||_{\infty}, ||f'||_{\infty} \leq 1$ from Lemma \ref{boundsol}
completes the proof. \hfill $\Box$

\vspace{.15in}

The quantity that is difficult to bound in practice when applying Theorem
\ref{atheorem1} is \[ \ee\left(D^2\ii[\abs{W-t}\leq\abs{D}]\right).\] One
tool which is useful for bounding this quantity is the following theorem.

\begin{theorem}\label{tool1} Assume that $W$ and $W'$ are real valued random
variables on the same probability space such that $\law(W')=\law(W)$. Let
$D= W'-W$. Then for any $t \in \rr$ and $c > 0$,
\[ \ee\bigl(D^2 \ii\{|W-t|\le |D|\}\bigr) \le 4 c \ee\bigl|\ee(D|W)\bigr| +
\ee\bigl(D^2 \ii\{|D|> c\}\bigr).
\]
\end{theorem}

However Theorem \ref{tool1} does not always give good bounds. The next
result, though more demanding, can lead to sharper bounds.

\begin{theorem}\label{tool2} Assume that $W$ and $W'$ are non-negative random
variables on the same probability space such that $\law(W')=\law(W)$; let
$D=W'-W$. Then for any positive constants $t$, $k_1$, $k_2$, $K_1$, $K_2$
and $K_3$ (where $k_2<k_1$ and $K_2<K_3$) we have
\begin{equation*}
\begin{split}
&\IE\left\{D^2 \ii[\abs{W-t}\leq \abs{D}]\right\}\leq
k_2+ k_1e_2 + e_1 + \frac{\IE\abs{\IE(D|W)}}{K_3-K_2} \times \\
&  \Bigg( k_2 \cdot \ln (k_1/k_2) + \sqrt{32tk_2k_1} + 2K_1t^{1/2}k_1+
4K_1\sqrt{k_1k_2}+ 4 K_1(tk_2k_1^3)^{1/4}\Bigg)
\end{split}
\end{equation*}
where
\begin{align*}
 e_1 &:= \IE\left\{\IE(D^2|W) \cdot
\II\left[\text{$\IE(D^2|W) > k_1$ \text{\textup{or}} $\IE(D^4|W) >
k_2(W+t)$}\right]
\right\}\\
 e_2 &:= \IP\left[\text{$\IE(D^2|W) < K_3$
    \text{\textup{or}} $\IE(D^4|W) > K_1^2 K_2 W$}\right]
\end{align*}
\end{theorem}

The following lemma will be used in the proofs of both Theorems
\ref{tool1} and \ref{tool2}.

\begin{lemma} \label{key}
Suppose that $W$ and $W'$ are random variables on the same probability
space such that $\law(W') = \law(W)$; set $D=W'-W$. Then, for any $a\leq b
\in \IR$ and $K>0$,
\begin{equation*}
    \ee\left(D^2 \ii[a\leq W\leq b,\abs{D}\leq K]\right) \leq
    (b-a+2K)\ee\abs{\ee(D|W)}
\end{equation*}
\end{lemma}

\begin{proof} Define
\begin{equation*}
 h(x) = \begin{cases}
         -\ahalf(b-a)-K & \text{if $x<a-K$,}\\
         x-\ahalf(a+b)  & \text{if $a-K\leq x\leq b+K$,}\\
         \ahalf(b-a)+K & \text{if $x>b+K$.}
        \end{cases}
\end{equation*}
and $H(x) := \int_0^x h(t)dt$. Observe that for any $0\leq s\leq 1$,
\begin{align*}
    \ii[a\leq W\leq b, \abs{D}\leq K]
        &\leq \ii[a-K\leq W+sD\leq b+K] \\
        &= h'(W+sD),
\end{align*}
and that
\begin{equation}    \label{bound1}
    \norm{h}_\infty = \ahalf(b-a)+K.
\end{equation}
Using Taylor expansion we have
\begin{equation} \begin{split} \label{lemmaequality1}
 0  & = \ee H(W') - \ee H(W)\\
    & = \ee(D h(W)) + \ee\left(D^2\int_0^1(1-s)h'(W+sD)ds\right),\\
\end{split} \end{equation} and thus
\begin{align*}
    &\ee\left(D^2 \ii[a\leq W\leq b,\abs{D}\leq K]\right)\\
    &\qquad = 2\ee\left(D^2\int_0^1(1-s)
            \ii[a\leq W\leq b,\abs{D}\leq K]ds\right)\\
    &\qquad \leq 2\ee\left(D^2\int_0^1(1-s)h'(W+sD)ds\right)\\
    &\qquad = -2\ee(D h(W))\kern10em\text{[by \eqref{lemmaequality1}]}\\
    &\qquad \leq 2\abs{\ee(\ee(D|W) h(W))}\\
    &\qquad \leq 2\norm{h}_\infty\ee\abs{\ee(D|W)}
\end{align*}
which together with \eqref{bound1} proves the claim.
\end{proof}

As the following argument shows, Theorem \ref{tool1} is a straightforward
consequence of Lemma \ref{key}.

{\sc Proof of Theorem \ref{tool1}} Clearly \begin{eqnarray*} \ee\bigl(D^2
\ii\{|W-t|\le |D|\}\bigr) & = & \ee\bigl(D^2 \ii\{|W-t|\le
|D|, |D|>c\}\bigr)\\
& & + \ \ee\bigl(D^2 \ii\{|W-t|\le |D|, |D| \leq c \}\bigr).
\end{eqnarray*} The first term is at most $\ee\bigl(D^2 \ii\{|D|>
c\}\bigr)$. To upper bound the second term, note that if $|W-t| \leq |D|
\leq c$, then $a \leq W \leq b$ where $a=t-c$ and $b=t+c$. Hence Lemma
\ref{key} gives that \[ \ee\bigl(D^2 \ii\{|W-t|\le |D|, |D| \leq c
\}\bigr) \leq 4 c \ee\bigl|\ee(D|W)\bigr|.\] \hfill $\Box$

We close this section by proving Theorem \ref{tool2}.

{\sc Proof of Theorem \ref{tool2}} Define
\begin{equation*}
    B(W) := \II\left[\IE [ (D^2 \II[\abs{D}> K_1 W^{1/2}] ) |W] \leq K_2,
            \IE (D^2|W) \geq K_3\right].
\end{equation*}
Now note that
\begin{equation*}
 \IE [(D^2\II[D^2> K_1^2 W])|W] \leq \frac{\IE (D^4|W)}{K_1^2 W}.
\end{equation*}
From this it is easy to see that
\begin{equation*}
 \IE(1-B(W)) \leq \IP\left[\text{$\IE (D^2|W) < K_3$
    or $\IE (D^4|W) > K_1^2 K_2 W$}\right] = e_2
\end{equation*}
Note that if $B(W)=1$ then
\begin{equation*}
\begin{split}
    \IE [(D^2 \II[\abs{D}\leq K_1 W^{1/2}])|W]
    & = \IE (D^2|W) - \IE [(D^2 \II[\abs{D}> K_1 W^{1/2}])|W]\\
    &\geq K_3-K_2.
\end{split}
\end{equation*}
Thus,
\begin{equation}
\begin{split}\label{concineq}
\IP[a\leq W\leq b]
&\leq \IE\left\{\II[a\leq W\leq b]B(W)\right\}  + e_2\\
&= \IE\left\{\frac{K_3-K_2}{K_3-K_2}\II[a\leq W\leq b]B(W)\right\}  + e_2\\
&\leq \frac{1}{K_3-K_2}\IE\left\{D^2\II[a\leq W\leq b,
             \abs{D}\leq K_1 W^{1/2}]\right\}+e_2\\
&\leq \frac{1}{K_3-K_2}\IE\left\{D^2\II[a\leq W\leq b,
             \abs{D}\leq K_1 b^{1/2}]\right\}+e_2\\
&\leq \frac{1}{K_3-K_2}(b-a+2K_1b^{1/2})\IE\abs{\IE(D|W)} + e_2
\end{split}
\end{equation}
where the last inequality is due to Lemma~\ref{key}.

Now, define
\begin{equation*}
 A(W) := \II\left[\IE(D^2|W) \leq k_1, \IE (D^4|W) \leq k_2(W+t)\right]
\end{equation*}
Then,
\begin{equation*}
\begin{split}
& \IE\{D^2(1-A(W))\} \\
&\qquad = \IE\left\{\IE (D^2|W) \cdot \II\left[\text{$\IE (D^2|W) > k_1$
or $\IE (D^4|W) > k_2(W+t)$}\right]
\right\}\\
&\qquad = e_1.
\end{split}
\end{equation*}
It follows that
\begin{equation*}
\begin{split}
&\IE\left\{D^2 \II[\abs{W-t}\leq \abs{D}]\right\} \\
&\qquad\leq \IE\left\{D^2 \II[(W-t)^2\leq D^2]A(W)\right\} + e_1\\
&\qquad\leq \IE\left\{\min\{D^2,D^4(W-t)^{-2}\} A(W)\right\} + e_1\\
&\qquad\leq \IE\left\{\min\{\IE (D^2|W),\IE(D^4|W) (W-t)^{-2}\} A(W)\right\} + e_1\\
&\qquad\leq \IE\left\{\min\{k_1,k_2(W+t)(W-t)^{-2}\}\right\} + e_1
\end{split}
\end{equation*}
Now,
\begin{equation*}
\begin{split}
&\IE\left\{\min\{k_1,k_2(W+t)(W-t)^{-2}\}\right\}\\
&\qquad = \int_0^\infty \IP[k_1\geq x, k_2(W+t)(W-t)^{-2}\geq x]dx\\
&\qquad = \int_0^{k_1} \IP[k_2(W+t)(W-t)^{-2}\geq x]dx \\
&\qquad \leq k_2 +  \int_{k_2}^{k_1} \IP[k_2(W+t)(W-t)^{-2}\geq x]dx
\end{split}
\end{equation*}
Suppose that $k_2(W+t)(W-t)^{-2}\geq x$. Then, solving the equation
$k_2(w+t)(w-t)^{-2} = x$, one has that
\begin{equation*}
\begin{split}
W & \in\left[t + \frac{k_2}{2x} - \sqrt{\frac{2tk_2}{x}+
\frac{k_2^2}{4x^2}},
 t + \frac{k_2}{2x} + \sqrt{\frac{2tk_2}{x}+ \frac{k_2^2}{4x^2}}\right]\\
& \subset\left[t - \sqrt{\frac{2tk_2}{x}}, t + \frac{k_2}{x} +
\sqrt{\frac{2tk_2}{x}}\right].
\end{split}
\end{equation*}
Thus, combining this with the concentration inequality \eqref{concineq},
\begin{equation*}
\begin{split}
&\IE\left\{D^2 \II[\abs{W-t}\leq \abs{D}]\right\}\\
&\leq k_2 + e_1 + \int_{k_2}^{k_1} \IP[k_2(W+t)(W-t)^{-2}\geq x]dx\\
&\leq k_2 +  e_1 + \int_{k_2}^{k_1} \IP\left[ t - \sqrt{\frac{2tk_2}{x}}
\leq W \leq
t + \frac{k_2}{x} + \sqrt{\frac{2tk_2}{x}}\right] dx\\
&\leq k_2 + k_1e_2 + e_1 + \frac{\IE\abs{\IE (D|W)
}}{K_3-K_2}\int_{k_2}^{k_1} \left(\frac{k_2}{x} +
\sqrt{\frac{8tk_2}{x}}+2K_1 \left(t + \frac{k_2}{x} +
\sqrt{\frac{2tk_2}{x}}\right)^{1/2}\right)
dx\\
&\leq k_2 + k_1e_2 + e_1 + \frac{\IE\abs{\IE (D|W)
}}{K_3-K_2}\int_{k_2}^{k_1} \left(\frac{k_2}{x} +
\sqrt{\frac{8tk_2}{x}}+2K_1t^{1/2} + 2K_1\sqrt{\frac{k_2}{x}} +
2K_1\left(\frac{2tk_2}{x}\right)^{1/4}\right)
dx\\
&\leq k_2+ k_1e_2 + e_1 + \frac{\IE\abs{\ee(D|W)}}{K_3-K_2} \times\\
&\Bigg( k_2(\ln k_1-\ln k_2) + \sqrt{32tk_2k_1} + 2K_1t^{1/2}k_1 +
4K_1\sqrt{k_1k_2}+ {\textstyle\frac{8}{3}} K_1(2tk_2)^{1/4}k_1^{3/4}
\Bigg).
\end{split}
\end{equation*}
This proves the claim. \hfill $\Box$

\section{Berry-Ess\'{e}en Bound for the Exponential Law: Version 2} \label{v2}

A main goal of this section is to prove Theorem \ref{basicexp} from the
introduction, and to develop tools for analyzing the error terms which
appear in it. In particular, the third term can be hard to bound. One way
to bound it is to apply Theorem \ref{tool1} from Section \ref{v1}. Another
way it to use the following more demanding result, which is analogous to
Theorem \ref{tool2} from Section \ref{v1}.

\begin{theorem}\label{tool3} Let $W$ and $W'$ be non-negative random variables
on the same probability space such that $\law(W')=\law(W)$. Suppose that
$\ee(D|W) = -\lambda (W-1)$, where $D = W'-W$ and $\lambda > 0$ is a fixed
constant. Then for any $t
> 0$ and $\kappa > 0$
\begin{align*}
\ee\bigl(D^2\ii\{|W-t|\le |D|\}\bigr) &\le 16\lambda^2\kappa^2 +
1040 \lambda^{3/2} \ee|W-1| \kappa \sqrt{t} \\
& \quad + 8\lambda \epsilon_2(\frc{1}{3}t) t + \epsilon_1(t)
\end{align*}
where $\epsilon_1$ and $\epsilon_2$ are functions defined on $(0,\infty)$
as
\begin{align*}
\epsilon_1(t) := \ee\bigl[\ee(D^2|W)\ii\{\ee(D^2|W) > 2\lambda(W + t)
\text{
 \textup{or} }
 \ee(D^4|W) > 4\lambda^2(\kappa^2W^2+\kappa^2t^2)\}\bigr]
\end{align*}
and
\[
\epsilon_2(t) := \pp\{\ee(D^2|W) < 2\lambda(W - \frc{1}{4}t) \text{
  \textup{or}
} \ee(D^4|W)  > 4\lambda^2(\kappa^2W^2+\kappa^2t^2)\}.
\] Moreover, the above bound holds if the assumption of positivity of $W$
is replaced by the assumption that $W$ is non-negative and assumes only
finitely many values.
\end{theorem}

{\it Remarks:} \begin{enumerate} \item The idea behind the formulation of
Theorem \ref{tool3} is the following: in many problems, we have
$\ee(D^4|W) \le 4\lambda^2(\kappa^2W^2 +\eta)$ where $\kappa$ is some
constant and $\eta$ is a negligible term (possibly random).

\item The random variable $W$ in the example of this paper can assume the value
0 with positive probability.
\end{enumerate}

\bigskip

It is easy to check by integration by parts that if a random variable $Z$
on $[0,\infty)$ is $\Exp(1)$, then $\ee[Zf'(Z)-(Z-1)f(Z)]=0$ for well
behaved functions $f$. This motivates the study of the solution $f(x)$ to
the equation \[ xf'(x)-(x-1)f(x) = \ii \{x \leq t \} - (1-e^{-t}) , \ \ x
\geq 0.\] Indeed, for such $f$ one has that
\[ \pp(W \leq t) - \pp(Z \leq t) = \ee[Wf'(W)-(W-1)f(W)], \] and the problem
becomes that of bounding \[ \ee[Wf'(W)-(W-1)f(W)]. \]

{\it Remark:} Earlier authors (Mann \cite{Man}, Luk \cite{Lu}, Pickett and
Reinert \cite{Re}) studied solutions of the equation
\[ xf''(x)-(x-1)f'(x) = h(x) - \int_0^{\infty} e^{-x} h(x), \] for
functions $h$ whose first $k$ derivatives are bounded. This is
complementary to our work, since our primary interest is in the function
$h(x) = \ii\{x \leq t\}$, which is not smooth.

\begin{lemma} \label{boundsol2} For every $t \in \rr$, the function
\[ f(x) := \frac{e^{-(t-x)^+} - e^{-t}}{x}, \ \ x > 0 \] satisfies
the equation
\begin{equation}\label{expeq}
x\fp(x) - (x-1)f(x) = \ii\{ x\le t\} - (1-e^{-t}), \ \ x\in \rr^+,
\end{equation} where $f'$ denotes the left-hand derivative of $f$.
Moreover, one has the bounds \[ ||f'||_{\infty} \leq t^{-1}, \ \
||f''||_{\infty} \leq max \{t^{-1},2t^{-2} \}.\] \end{lemma}

\begin{proof} Clearly, $f$ is infinitely differentiable on $\rr^+ \backslash \{t\}$.
The left-hand and right-hand derivatives at $t$ exist and are unequal,
which is why we let $f'$ denote the left-hand derivative of $f$. Then for
$0< x \le t$,
\begin{equation}\label{fp1}
\begin{split}
\fp(x) &= \frac{d}{dx} \biggl(\frac{e^{-(t-x)} - e^{-t}}{x}\biggr) =
\frac{xe^x - e^x + 1}{x^2} e^{-t}
\end{split}
\end{equation}
which gives
\begin{align*}
x\fp(x) - (x-1)f(x) &= 1 - (1-e^{-t}).
\end{align*}
Similarly, for $x>t$,
\begin{equation}\label{fp2}
\begin{split}
\fp(x) &= \frac{d}{dx}\biggl(\frac{1-e^{-t}}{x}\biggr) =
\frac{e^{-t}-1}{x^2}
\end{split}
\end{equation}
which gives
\begin{align*}
x\fp(x) - (x-1)f(x) &= - (1-e^{-t}).
\end{align*}
Thus, the function $f$ is a solution to (\ref{expeq}).

The easiest way to get a uniform bound on $\fp$ is perhaps by directly
expanding in power series. When $0< x \le t$, we recall (\ref{fp1}) to get
\begin{equation}\label{fp3}
\fp(x) = \frac{xe^x - e^x + 1}{x^2} e^{-t} = e^{-t} \sum_{k=0}^\infty
\frac{x^k}{k! (k+2)}.
\end{equation}
This shows that for $x \in (0,t]$,
\[
0 \le \fp(x) \le e^{-t}\sum_{k=0}^\infty \frac{t^k}{(k+1)!} =
\frac{1-e^{-t}}{t} \le \min\{1,t^{-1}\}.
\]
Again, for $x > t$, we directly see from (\ref{fp2}) that $\fp(x) \le 0$
and
\begin{equation*}
|\fp(x)| \le \frac{1-e^{-t}}{t^2}
\end{equation*}
Combining, we get
\begin{equation*}
\|\fp\|_\infty \le t^{-1}.
\end{equation*}

Now, $\fp$ is positive in $(0,t]$ and negative in $(t,\infty)$. Therefore
$f$ attains its maximum at $t$. It is now easy to see that for all $x>0$,
\begin{equation*}
0\le f(x) \le \frac{1-e^{-t}}{t} \le 1.
\end{equation*}
Using (\ref{fp3}) we see that for $0<x \le t$
\begin{equation*}
\begin{split}
0\le \fpp(x) &= e^{-t} \sum_{k=1}^\infty \frac{kx^{k-1}}{k!(k+2)} = e^{-t}
\sum_{k=0}^\infty \frac{x^k}{k!(k+3)} \le \frac{1-e^{-t}}{t}
\end{split}
\end{equation*}
and for $x> t$,
\begin{equation*}
0 \le \fpp(x) = \frac{2(1-e^{-t})}{x^3} \le \frac{2(1-e^{-t})}{t^3}.
\end{equation*}
Combining, we get, for all $x > 0$,
\begin{equation*}
0 \le \fpp(x) \le \max\{t^{-1},2t^{-2}\}.
\end{equation*} This completes the proof. \end{proof}

Now the main results of this section will be proved.

{\sc Proof of Theorem \ref{basicexp} } Fix $t>0$ and consider the Stein
equation
\begin{equation}\label{a}
    x f'(x) - (x-1)f(x) = \ii[x\leq t] - \IP[Z\leq t]
\end{equation}
for $x>0$, where $Z\sim\Exp(1)$. From Lemma \ref{boundsol2}, its solution
$f$ satisfies the non-uniform bounds
\begin{equation} \label{bds}
    \norm{f'}_{\infty} \leq t^{-1},\ \norm{f''}_{\infty} \leq \max\{t^{-1},2t^{-2}\},
\end{equation}
where $f''$ denotes the left derivative of $f'$, as $f'$ has a
discontinuity in $t$.

Assume first that $W$ is positive. Defining $G(w)=\int_0^w f(x) dx$,
Taylor's expansion gives that
\begin{equation*}
    G(W') = G(W) + D f(W) + D^2\int_0^1 (1-s) f'(W+sD)\,ds.
\end{equation*} The hypothesis $\ee(D|W)=-\lambda(W-1)$ gives that
\begin{equation*}
\begin{split}
    0  & = \IE\left\{G(W') - G(W)\right\}\\
     & = \IE\left\{-\lambda(W-1) f(W)\right\}
    + \IE\left\{D^2\int_0^1 (1-s) f'(W+s D)\,d s\right\},
\end{split}
\end{equation*}
and hence
\begin{equation*}
    \IE\left\{(W-1) f(W)\right\}
    = \IE\left\{\frac{1}{\lambda}D^2\int_0^1 (1-s) f'(W+sD)\,d s\right\}.
\end{equation*}

Taking expectation on (\ref{a}) with respect to $W$, we thus have
\begin{equation} \label{b}
\begin{split}
    &\IP[W\leq t] - \IP[Z\leq t]\\
    &\qquad = \IE\left\{W f'(W) - (W-1)f(W)\right\} \\
    &\qquad = \IE\left\{W f'(W)
        - \frac{1}{\lambda}D^2\int_0^1 (1-s) f'(W+s D)\,d s\right\}\\
    &\qquad = \IE\left\{\left(W - D^2/(2\lambda)\right) f'(W)\right\}\\
    &\qquad\qquad\qquad + \IE\left\{\frac{ D^2}{\lambda}
       \left[ \frac{1}{2} f'(W)
        - \int_0^1 (1-s) f'(W+s D)\,ds \right] \right\}\\
    &\qquad = \IE\left\{\left(W - D^2/(2\lambda)\right) f'(W)\right\}\\
    &\qquad\qquad\qquad + \IE\left\{\frac{ D^2}{\lambda}
        \int_0^1 (1-s) (f'(W)-f'(W+s D)) d s \right\}.
\end{split}
\end{equation}
Note now that for any $x,y>0$,
\begin{equation*}
    \abs{f'(x)-f'(y)} \leq
    \begin{cases}
    \norm{f''}_{\infty}|x-y| &
        \text{if $x$ and $y$ lie on the same side of $t$,}\\
    2\norm{f'}_{\infty} &\text{otherwise.}
    \end{cases}
\end{equation*}
Also, if $x$ and $y$ lie on different sides of $t$, then $|x-t|\leq|x-y|$.
Thus
\begin{equation*}
\begin{split}
    &\int_0^1 (1-s)\left|f'(W)-f'(W+s D)\right| d s \\
    &\quad \leq \norm{f''}_{\infty} \int_0^1 (1-s)s|D|\, d s
    + 2\norm{f'}_{\infty} \int_0^1 (1-s) \ii \left[\abs{W-t}\leq \abs{s D}\right] d s\\
    &\quad \leq \frac{1}{6}|D|\norm{f''}_{\infty}
    + \norm{f'}_{\infty} \ii \left[ \abs{W-t}\leq |D|\right].
\end{split}
\end{equation*}
Putting the steps together we obtain from (\ref{b})
\begin{eqnarray*}
    \left|\IP[W\leq t] - \IP[Z\leq t]\right| & \leq &
    \norm{f'}_{\infty} \IE\left|W-\frac{\ee(D^2|W)}{2\lambda} \right| \\
    & & + \frac{1}{6\lambda}\norm{f''}_{\infty} \IE|D|^3
    + \frac{1}{\lambda}\norm{f'}_{\infty} \IE\left\{D^2 \ii[\abs{W-t}
    \leq \abs{D}]\right\},
\end{eqnarray*}
and with the bounds (\ref{bds}) the claim follows for positive $W$.

    To treat the case where $W$ can also equal 0, choose
    $0<\delta<1$ and define $W_{\delta}:=(1-\delta)W+\delta$,
    $W'_{\delta}:=(1-\delta)W'+\delta$ and
    $t_{\delta}:=(1-\delta)t+\delta$. One sees that $W_{\delta}$
    is a positive random variable, and that
    $\ee(D_{\delta}|W_{\delta}) = -\lambda (W_{\delta}-1)$ where
    $\lambda$ is the same as for the pair $(W,W')$. Moreover $\pp
    \{W \leq t \} = \pp \{W_{\delta} \leq t_{\delta} \}$, so it
    follows that \begin{equation*}
\begin{split}
|\pp\{W\le t\} - \pp\{Z\le t\}| &\le \frac{1}{2\lambda t_{\delta}} \ee
\bigl|2\lambda W_{\delta} -\ee(D_{\delta}^2|W)\bigr| +
 \frac{\max\{t_{\delta}^{-1},2t_{\delta}^{-2}\}}{4\lambda}
\ee|D_{\delta}|^3 \\
& \quad + \frac{1}{\lambda t_{\delta}} \ee \bigl(D_{\delta}^2
\ii\{|W_{\delta}-t_{\delta}|\le |D_{\delta}|\}\bigr).
\end{split}
\end{equation*}Since $D_{\delta}=(1-\delta)D$, the first two error terms
 are continuous in $\delta$ and converge
to the corresponding error terms for $W$
 when $\delta \rightarrow 0$. The same is true for the third error term,
 as can be seen from the fact that $|W_{\delta}-t_{\delta}| \leq |D_{\delta}|$
if and only if $|W-t| \leq |D|$. This completes the proof. \hfill $\Box$ \\

Next, we prove Theorem \ref{tool3}.

{\sc Proof of Theorem \ref{tool3} } First we treat the case that $W$ is
always positive. Throughout we shall be using $V:= (2
\lambda)^{-1/2}(W'-W)$ instead of $D (=W'-W)$, simply because $D$ occurs
with a factor of $(2 \lambda)^{-1/2}$ attached with it on most occasions.

Suppose for each $0<s\le t$, we have numbers $u(s,t)$ and $v(s,t)$ such
that whenever $s\le a\le b \le t$, we have
\[
\pp\{a\le W\le b\} \le u(s,t)(b-a) + v(s,t).
\]
Fix $t\in \rr$. Let $A(W)= \ii\{\ee(V^2|W) \le W+t, \ \ee(V^4|W) \le
\kappa^2 (W+ t)^2\}$. Then
\begin{equation*}
\begin{split}
& \ee\bigl(V^2(1-A(W))\bigr)\\
& \leq \ee\bigl(\ee(V^2|W)\ii\{\ee(V^2|W) > W + t \text{ or } \ee(V^4|W)
> \kappa^2W^2+\kappa^2t^2\}\bigr) \\
&=: e_1(t).
\end{split}
\end{equation*}
It follows that
\begin{align*}
&\ee(V^2 \ii\{|W-t| \le |W'-W|\}) \\
& \leq \ee(V^2 \ii\{(W'-W)^2 \ge (W-t)^2\} A(W)) + e_1(t)
\\
&\le \ee(\min\{2\lambda (W-t)^{-2}V^4, V^2\}A(W)) + e_1(t)
\\
&\le \ee(\min\{2\lambda (W-t)^{-2}\ee(V^4|W), \ee(V^2|W)\}A(W)) + e_1(t)
\\
&\le \ee(\min\{2\lambda\kappa^2(W-t)^{-2}(W+t)^2, W+t\}) + e_1(t).
\end{align*}
Now
\begin{equation}\label{int1}
\begin{split}
& \ee(\min\{2\lambda\kappa^2(W-t)^{-2}(W+t)^2,W+t\})\\
&=\int_0^\infty \pp\{2\lambda \kappa^2(W-t)^{-2}(W+t)^2 \ge x,  W + t \ge
x\}\;dx.
\end{split}
\end{equation}
Now take any $x \ge 8\lambda\kappa^2$. Let $c(x) =
\sqrt{\frac{2\lambda\kappa^2}{x}}$. Then the following are easily seen to
be equivalent:
\begin{align*}
2\lambda \kappa^2(W-t)^{-2}(W+t)^2  \ge x &\iff |W-t| \le
c(x)(W+t) \\
&\iff \frac{1-c(x)}{1+c(x)} t \le W\le \frac{1+c(x)}{1-c(x)}t.
\end{align*}
Let $a(x)= (1-c(x))/(1+c(x))$ and $b(x)=(1+c(x))/(1-c(x))$. Note that
since $x\ge 8\lambda\kappa^2$, therefore $c(x) \le 1/2$ and so $a(x) \ge
1/3$, $b(x) \le 3$, and
\[
b(x) - a(x) = \frac{4c(x)}{1-c(x)^2} \le \frc{16}{3}c(x).
\]
Now if $W\le 3t$ then $W+t\le 4t$. Thus, the integrand in (\ref{int1}) is
zero for $x > 4t$. Combining, we see that
\begin{equation}\label{term1}
\begin{split}
& \ee(\min\{2\lambda(W-t)^{-2}(W+t)^2,W+t\}) \\
&\le 8\lambda\kappa^2 + \int_{8\lambda\kappa^2}^{4t} \pp\{a(x)t \le
W\le b(x)t\}dx \\
&\le 8\lambda\kappa^2 + \int_{8 \lambda\kappa^2}^{4t} \bigl(u(\frc{1}{3}t,
3t)\frc{16}{3}c(x)t + v(\frc{1}{3}t,
3t)\bigr)dx\\
&\le 8\lambda\kappa^2 + 22 u(\frc{1}{3}t, 3t) t\kappa \sqrt{2\lambda t} +
4tv(\frc{1}{3}t, 3t).
\end{split}
\end{equation}

Next, we proceed to find suitable values of $u(s,t)$ and $v(s,t)$. Fix $0
< s \le a \le b \le t$. Let
\[
B(W) = \ii\{\ee(V^2\ii\{|V|> 2\kappa\sqrt{W+a}\}|W) \le \frc{1}{4}(W+a), \
\ee(V^2|W) \ge W-\frc{1}{4}a \}.
\]
Now note that
\[
\ee(V^2\ii\{|V|> 2\kappa\sqrt{W+a}\}|W) \le
\frac{\ee(V^4|W)}{4\kappa^2(W+a)}.
\]
From this it is easy to see that
\begin{equation*}
\begin{split}
\ee(1-B(W)) &\le \pp\{\ee(V^2|W) < W - \frc{1}{4}a \text{ or } \ee(V^4|W)
> \kappa^2W^2+\kappa^2a^2\} \\
&=: e_2(a)
\end{split}
\end{equation*}
Note that if $B(W)=1$ then $\ee(V^2\ii\{|V|\le 2 \kappa\sqrt{W+a}\}|W)\ge
\frac{3}{4}W - \frac{1}{2}a$. So, if $W\ge a$ and $B(W) = 1$,
$\ee(V^2\ii\{|V|\le 2 \kappa\sqrt{W+a}\}|W)\ge \frac{1}{4}a$. Thus,
\begin{align*}
a\pp\{a \le W \le b\} &\le 4\ee\bigl(V^2\ii\{a\le W\le b, \ |V| \le 2
\kappa\sqrt{W+a}\}\bigr)
+ ae_2(a) \\
&\le  4\ee\bigl(V^2\ii\{a\le W\le b, \ |V| \le 2 \kappa\sqrt{b+a}\}\bigr)
+ a e_2(a) \\
&= 2\lambda^{-1} \ee\bigl(D^2 \ii\{a\le W\le b, \ |D| \le 2 \kappa\sqrt{2
\lambda (b+a)}\}\bigr)  +  ae_2(a).
\end{align*}
where $D = W'-W$. Using Lemma \ref{key}, we get
\[
a\pp\{a \le W \le b\} \le 2\bigl(b-a + 4 \kappa\sqrt{2 \lambda
(b+a)}\bigr) \ee|W-1| +  a e_2(a).
\]
Finally, note that $e_2$ is a monotonically decreasing function. Thus, we
can take
\[
u(s,t) = \frac{2 \ee|W-1|}{s}
\]
and
\[
v(s,t) =\frac{16 \kappa\sqrt{\lambda t}\ee|W-1|}{s} + e_2(s).
\]
Using these expressions for $u$ and $v$ in (\ref{term1}), we get
\[
\ee\bigl(V^2\ii\{|W-t|\le |W'-W|\}\bigr) \le 8\lambda\kappa^2 + 520
\ee|W-1| \kappa \sqrt{\lambda t} + 4e_2(\frc{1}{3}t) t + e_1(t).
\]
Put $\epsilon_1(t) = 2\lambda e_1(t)$ and $\epsilon_2(t) = e_2(t)$ to get
the final expression in Theorem \ref{tool3}.

    Finally, suppose that $W$ might take the value 0, but that $W$
    assumes only finitely many values. As in the proof of Theorem
    \ref{basicexp}, for $0<\delta<1$ define
    $W_{\delta}:=(1-\delta)W+\delta$,
    $W'_{\delta}:=(1-\delta)W'+\delta$ and
    $t_{\delta}:=(1-\delta)t+\delta$. Since
    $|W_{\delta}-t_{\delta}| \leq |D_{\delta}|$ if and only if
    $|W-t| \leq |D|$, it follows that $\ee(D^2 \ii \{|W-t| \leq
    |D| \})$ is the limit as $\delta \rightarrow 0$ of
    $\ee(D_{\delta}^2 \ii \{|W_{\delta}-t_{\delta}| \leq
    |D_{\delta}| \})$. It is easily checked that
    $\ee(D_{\delta}^2|W)>2 \lambda(W_{\delta}+t_{\delta})$ implies
    that $\ee(D^2|W)>2 \lambda(W+t)$ and that $\ee(D_{\delta}^4|W)
    > 4 \lambda^2(\kappa^2 W_{\delta}^2+\kappa^2t_{\delta}^2)$
    implies that $\ee(D^4|W) > 4 \lambda^2(\kappa^2
    W^2+\kappa^2t^2)$. We claim that $\ee(D_{\delta}^2|W)<2
    \lambda(W_{\delta}-\frac{1}{4}t_{\delta})$ implies that
    $\ee(D^2|W)<2 \lambda(W-\frac{1}{4}t)$ provided that $\delta$
    is sufficiently small. Indeed, since $W$ takes only finitely
    many values, there is an $m_t>0$ such that $\ee(D^2|W)<2
    \lambda(W-\frac{1}{4}t)$ if and only if $\ee(D^2|W)<2
    \lambda(W-\frac{1}{4}t)+m_t$. The claim now follows since
    $\ee(D_{\delta}^2|W)<2
    \lambda(W_{\delta}-\frac{1}{4}t_{\delta})$ implies that
    $\ee(D^2|W)<2 \lambda(W-\frac{1}{4}t) + \frac{3 \lambda
    \delta}{2(1-\delta)} + \delta \ee(D^2|W)$. Hence the theorem
    follows by letting $\delta \rightarrow 0$.  \hfill $\Box$

\section{Example: Spectrum of Bernoulli-Laplace chain} \label{n/2}

This section proves Theorem \ref{sharper} of the introduction. Throughout
we let $W$ denote the random variable defined by \[ W(i) :=
\frac{(n-2i)(n+2-2i)}{2n}, \] where $n$ is even and $i \in
\{0,1,\cdots,\frac{n}{2} \}$ is chosen with probability $\pi(i)$ equal to
\[  \frac{{n \choose i} - {n \choose i-1}}{{n \choose n/2}} \ \ \ \mbox{if
\ $1 \leq i \leq \frac{n}{2}$}, \ \ \  \frac{1}{{n \choose n/2}} \ \ \
\mbox{if $i=0$.} \]

Letting $Z$ be an $\Exp(1)$ random variable and $C$ a universal constant,
the upper bound
\[ |\pp(W \leq t) - \pp(Z \leq t)| \leq \frac{C}{\sqrt{n}} \] will be
proved in two steps. Subsection \ref{smallt} uses the machinery of Section
\ref{v1} to treat the case that $t \leq 1$, and Subsection \ref{larget}
uses the machinery of Section \ref{v2} to treat the case that $t \geq 1$.
One interesting feature of the proof is that the exchangeable pairs used
in these two subsections are different (but closely related). We also show
(in Subsection \ref{smallt}), that combining the machinery of Section
\ref{v1} with a concentration inequality, one can obtain, with less
effort, a slightly weaker $O(\frac{\log(n)}{\sqrt{n}})$ upper bound.

Finally, Subsection \ref{lower} shows that the $O(n^{-1/2})$ rate is
sharp, by constructing a sequence of $n$'s tending to infinity and
corresponding $t_n$'s such that \[ |\pp(W_n \leq t_n) - \pp(Z \leq t_n)| =
\frac{2e^{-2}}{\sqrt{n}} + O(1/n).\]

\subsection{Upper bound for small $t$} \label{smallt}

The purpose of this subsection is to use the machinery of Section \ref{v1}
to prove Proposition \ref{small}, which implies the upper bound of Theorem
\ref{sharper} of the introduction for $t \leq 1$.

\begin{prop} \label{small}
\[ |\pp(W \leq t) - \pp(Z \leq t)| \leq \frac{ C \cdot
\max(1,t^{1/2})}{\sqrt{n}} \] for a universal constant $C$.
\end{prop}

To begin we define an exchangeable pair $(W,W')$ and perform some
computations with it. The definition of $(W,W')$ and the fact that the
computations work out so neatly may seem unmotivated. There is an
algebraic motivation for our choices, and so as not to interrupt our
self-contained probabilistic treatment, we explain this in the appendix.

To construct an exchangeable pair $(W,W')$, we specify a Markov chain $K$
on the set $\{0,1,\cdots,\frac{n}{2} \}$ which is reversible with respect
to $\pi$. This means that $\pi(i)K(i,j)= \pi(j)K(j,i)$ for all $i,j$.
Given such a Markov chain $K$, one obtains the pair $(W,W')$ in the usual
way (see for instance \cite{RR}): choose $i$ from $\pi$, let $W=W(i)$, and
let $W'=W(j)$, where $j$ is obtained from $i$ by taking one step using the
Markov chain $K$.

The Markov chain which turns out to be useful is a birth-death chain on
$\{0,1,\cdots,\frac{n}{2} \}$ where the transition probabilities are
\[ K(i,i+1):= \frac{n-i+1}{n(n-2i)(n-2i+1)} \]
\[ K(i,i-1):=\frac{i}{n(n-2i+1)(n-2i+2)} \]
\[ K(i,i):=1-K(i,i+1)-K(i,i-1), \] with the exception of $K(i,i+1)$ if
$i=n/2$, which we define to be zero.

It is easily checked that $K$ is reversible with respect to $\pi$, so the
resulting pair $(W,W')$ is exchangeable. (In fact the machinery of Section
\ref{v1} only uses that $W$ and $W'$ have the same law, which follows from
the fact that $K$ has $\pi$ as a stationary distribution, but the
exchangeability is good to record).

Lemma \ref{momcom} performs some moment computations related to the pair
$(W,W')$.

\begin{lemma} \label{momcom} Letting $D:=W'-W$, one has that:
\begin{enumerate}
\item $\ee(D|W)=-\frac{2}{n^2}$ if $W \neq 0$; $\ee(D|W=0)=\frac{1}{n}$.
\item $\ee(W)=1$.
\item $\ee(D^2|W)=\frac{4}{n^2}$.
\item $\ee(D^4|W)=\left( \frac{32}{n^3}-\frac{64}{n^4} \right)W +
\frac{64}{n^4}$.
\item $\ee(D^4)=\frac{32}{n^3}$.
\end{enumerate}
\end{lemma}

\begin{proof} Since $i$ is determined by $W(i)$, conditional expectations given
$W$ can be computed using conditional expectations given $i$. Supposing
that $i \neq n/2$, \begin{eqnarray*} \ee(D|i) & = & K(i,i+1)(W(i+1)-W(i))
+ K(i,i-1)(W(i-1)-W(i)) \\
& = & \frac{n-i+1}{n(n-2i)(n-2i+1)} \frac{2(2i-n)}{n} +
\frac{i}{n(n-2i+1)(n-2i+2)} \frac{2(n-2i+2)}{n} \\
& = & - \frac{2}{n^2}. \end{eqnarray*} If $i=n/2$, then $\ee(D|i) =
K(i,i-1) (W(i-1)-W(i)) = \frac{1}{n}$, so part 1 is proved.

For part 2, argue as in part 1 (separately treating the cases $i \neq n/2$
and $i=n/2$) to compute that $\ee(D^3|W)= -\frac{16}{n^3}(W-1)$. Since $W$
and $W'$ are exchangeable, $\ee(D^3)=0$. Thus \[ \ee(W-1)= -\frac{n^3}{16}
\ee[\ee(D^3|W)] = -\frac{n^3}{16} \ee(D^3) = 0,\] so $\ee(W)=1$.

For parts 3 and 4, one argues as in part 1 to compute both sides
(separately treating the cases $i \neq n/2$ and $i=n/2$) and checks that
they are equal. For part 5, note that \[ \ee(D^4)=\ee[\ee(D^4|W)]= \left(
\frac{32}{n^3} - \frac{64}{n^4} \right) \ee[W] + \frac{64}{n^4} =
\frac{32}{n^3}, \] where the final equality is part 2. \end{proof}

Using these moment computations, we deduce Proposition \ref{small}.

{\sc Proof of Proposition \ref{small}} We apply Theorem \ref{atheorem1} to
the pair $(W,W')$ with the value $\lambda=\frac{2}{n^2}$. Then the first
two error terms actually vanish. Indeed, part 1 of Lemma \ref{momcom}
gives that
\[ \IE\abs{(\lambda^{-1}\IE(D|W)+1)\ii[W>0]} = 0, \] and part 3 of Lemma
\ref{momcom} gives that \[ \IE\abs{{\textstyle\frac{1}{2\lambda}}
\IE(D^2|W)-1} = 0.\]

To analyze the third error term, use the Cauchy-Schwarz inequality and
parts 3 and 4 of Lemma \ref{momcom} to obtain that
\[ \frac{n^2}{12} \ee|D^3| \leq \frac{n^2}{12} \sqrt{\ee(D^2) \ee(D^4)} =
\sqrt{\frac{8}{9n}}.\]

To bound the fourth error term, apply Theorem \ref{tool2} with \[
k_1=\frac{4}{n^2}, \ k_2=\frac{48}{n^3}, \ K_1=\sqrt{\frac{48}{n}}, \
K_2=\frac{1}{n^2}, \ K_3=\frac{4}{n^2}.\] Note (as required by the
theorem), that $K_2<K_3$ and that for $n>12$, $k_2<k_1$. From part 1 of
Lemma \ref{momcom} and the fact that $\pp(W=0)=\frac{2}{n+2}$, one
computes that $\ee|\ee(D|W)|=\frac{4}{n(n+2)}$.

It is necessary to upper bound \[ e_1 = \IE\left\{\IE(D^2|W) \cdot
\II\left[\text{$\IE(D^2|W)
> k_1$ or $\IE(D^4|W) > k_2(W+t)$}\right] \right\}. \] Note from part 3
of Lemma \ref{momcom} that \[ \pp \left[ \ee(D^2|W)>\frac{4}{n^2} \right]
=0
\] and from part 4 of Lemma \ref{momcom} that
\begin{equation} \begin{split} \label{bound} \pp \left[ \ee(D^4|W)>\frac{48}{n^3}(W+t)
\right] & \leq  \pp \left[
\ee(D^4|W)>\frac{48}{n^3} W \right] \\
& =  \pp(W<4/(n+4)) \\
& =  \pp(W=0) \\
& =  \frac{2}{n+2}.\\ \end{split} \end{equation} Thus \[ e_1 \leq
\frac{4}{n^2} \pp(W=0) = \frac{8}{n^2(n+2)}.\] It is also necessary to
upper bound \[
 e_2 = \IP\left[\text{$\IE(D^2|W) < K_3$
    or $\IE(D^4|W) > K_1^2 K_2 W$}\right] .\] Note from part 3
of Lemma \ref{momcom} that \[ \pp \left[ \ee(D^2|W) < \frac{4}{n^2}
\right] =0
\] and from (\ref{bound}) that $\pp[\ee(D^4|W)>\frac{48}{n^3}
W]=\frac{2}{n+2}$. Thus $e_2 = \frac{2}{n+2}$. Plugging into Theorem
\ref{tool2}, one obtains that \[ \frac{n^2}{4} \ee\bigl(D^2 \ii\{|W-t|\le
|D|\}\bigr) \leq \frac{C \cdot \max \{1,t^{1/2}\}}{\sqrt{n}},\] for a
universal constant $C$. This completes the proof. \hfill $\Box$

\vspace{.15in}

To close this subsection, we show how the machinery of Section \ref{v1},
together with a concentration inequality for $W'-W$, leads to a simpler
proof (avoiding the use of Theorem \ref{tool2}) that
\[ |\pp(W \leq t) - \pp(Z \leq t)| \leq C \sqrt{\frac{\log(n)}{n}},
\] for a universal constant $C$. We hope that this approach will be useful
in other settings (a concentration inequality for $W'-W$ can be very
useful for normal approximation by Stein's method; see the survey
\cite{CS}).

The following lemma is helpful for obtaining a concentration result for
$W'-W$.

\begin{lemma}
 \label{binomial} Let $a$ be an integer such that $0 \leq a \leq
 \frac{n}{2}$. Then ${n \choose \frac{n}{2}-a}/{n \choose
 \frac{n}{2}} \leq e^{-\frac{a(a-1)}{n}}$.
\end{lemma}

\begin{proof} The  result is visibly true for $a = 0$, so suppose that $a \geq
1$. Observe that \begin{eqnarray*}  \frac{{n \choose \frac{n}{2}-a}}{{n \choose
 \frac{n}{2}}}  & = & \frac{(\frac{n}{2}) \cdots
(\frac{n}{2}-a+1)}{(\frac{n}{2}+1) \cdots (\frac{n}{2}+a)}\\ & \leq &
\frac{(\frac{n}{2}) \cdots (\frac{n}{2}-a+1)}{(\frac{n}{2})^a}\\ & = &
\prod_{i=1}^{a-1} (1-\frac{2i}{n})\\
& = & e^{\sum_{i=1}^{a-1} \log(1-\frac{2i}{n}) }\\ & \leq & e^{-
\sum_{i=1}^{a-1} \frac{2i}{n} }\\
& = & e^{-\frac{a(a-1)}{n}}. \end{eqnarray*} \end{proof}

Proposition \ref{concen} gives the concentration inequality for $W'-W$. As
usual $\lceil x \rceil$ denotes the smallest integer greater than or equal
to $x$.

\begin{prop} \label{concen} $\pp(|W'-W| > c) \leq n^{-5/2}$
for $c= \frac{4}{n} \left( \left \lceil \sqrt{ \frac{5}{2} n \log(n)}
\right \rceil + 1 \right)$.
\end{prop}

\begin{proof} Since the Markov chain $K$ used to construct $(W,W')$ is a birth death
chain, it is easily checked from the definition of $W$ that $|W'(i)-W(i)|
\leq \frac{2}{n}(n-2i+2)$ for all $i$. Thus for $c$ as in the proposition,
\begin{eqnarray*}
\pp(|W'-W| > c ) & \leq & \pp \left( \frac{2}{n}(n-2i+2)
> c \right)\\
& = & \pp \left( i<\frac{n}{2}-\frac{cn}{4} +1 \right)\\
& = & \pp \left[ i<\frac{n}{2}+1-\left( \left \lceil \sqrt{ \frac{5}{2} n
\log(n)} \right \rceil + 1 \right) \right] .
\end{eqnarray*} From the definition of the probability measure $\pi$, it
is clear that for integral $a$, $\pp(i<\frac{n}{2}+1-a) = \frac{{n \choose
\frac{n}{2}-a}}{{n \choose \frac{n}{2}}}$. Hence the proposition follows
from Lemma \ref{binomial}.
\end{proof}

This leads to the following proposition.

\begin{prop}
\[ |\pp(W \leq t) - \pp(Z \leq t)| \leq C \sqrt{\frac{\log(n)}{n}},
\] for a universal constant $C$. \end{prop}

\begin{proof} As in the proof of Proposition \ref{small}, apply Theorem
\ref{atheorem1} to the pair $(W,W')$ with the value
$\lambda=\frac{2}{n^2}$. The first three terms are bounded as in the proof
of Proposition \ref{small}. To bound the fourth term, note from Theorem
\ref{tool1} that
\[ \frac{1}{2 \lambda} \ee\bigl(D^2 \ii\{|W-t|\le |D|\}\bigr) \le
n^2 c \ee\bigl|\ee(D|W)\bigr| + \frac{n^2}{4} \ee\bigl(D^2 \ii\{|D|>
c\}\bigr)
\] for any $c>0$. From part 1 of Lemma \ref{momcom} one computes that
$\ee|\ee(D|W)|=\frac{4}{n(n+2)}$. One checks from the definitions that
$|W'-W| \leq 2+\frac{4}{n}$, so that $(W'-W)^2 \leq 16$ since $n$ is even.
Choosing $c=\frac{4}{n} \left( \left \lceil \sqrt{ \frac{5}{2} n \log(n)}
\right \rceil + 1 \right)$, it follows from Proposition \ref{concen} that
\[ \ee\bigl(D^2 \ii\{|D|> c\}\bigr) \leq 16 \pp(|D|>c) \leq 16 n^{-5/2}.\]
This proves the proposition. \end{proof}

\subsection{Upper bound for large $t$} \label{larget}

The purpose of this subsection is to apply the machinery of Section
\ref{v2} to prove the following Proposition, which gives the upper bound
in Theorem \ref{sharper} in the introduction for $t \geq 1$.

\begin{prop} \label{large}
\[ |\pp(W \leq t) - \pp(Z \leq t)| \leq \frac{ C \cdot
\max(1,t^{-3})}{\sqrt{n}} \] for a universal constant $C$.
\end{prop}

The pair $(W,W')$ used in this subsection is somewhat different from the
pair used in Subsection \ref{smallt}; for a discussion of the relationship
between the two pairs see the remark below. As with the pair from
Subsection \ref{smallt}, the definition and the fact that the computations
work out so nicely may seem unmotivated. The algebraic motivation for the
choices is discussed in the appendix.

To construct an exchangeable pair $(W,W')$, we specify a Markov chain $K$
on the set $\{0,1,\cdots,\frac{n}{2} \}$ which is reversible with respect
to $\pi$ (i.e. one has that $\pi(i)K(i,j)= \pi(j)K(j,i)$ for all $i,j$).
Given such a Markov chain $K$, one obtains the pair $(W,W')$ by choosing
$i$ from $\pi$, letting $W=W(i)$, and setting $W'=W(j)$, where $j$ is
obtained from $i$ by taking one step using the Markov chain $K$.

The Markov chain which turns out to be useful is a birth-death chain on
$\{0,1,\cdots,\frac{n}{2} \}$ whose only non-zero transition probabilities
are \[ K(i,i+1):= \frac{(n-i+1)(n-2i)}{n(n-2i+1)}, \
K(i,i-1):=\frac{i(n-2i+2)} {n(n-2i+1)}. \] It is easily checked that $K$
is reversible with respect to $\pi$ so that $(W,W')$ is exchangeable. (In
fact the machinery of Section \ref{v2} only uses that $W$ and $W'$ have
the same law).

{\it Remark:} If $K(i,j)$ denotes the transition probabilities of this
subsection, and $\tilde{K}(i,j)$ denotes the transition probabilities from
Subsection \ref{smallt}, one can verify the relation
\[ \tilde{K}(i,j) = \frac{4}{n^2} \frac{K(i,j)}{(W(i)-W(j))^2} , \ \
\forall i \neq j.\] Letting $D=W'-W$ for the pair of this subsection and
$\tilde{D},\tilde{W}$ the corresponding quantities for the pair from
Subsection \ref{smallt}, it follows that \[
\ee[\tilde{D}^r|\tilde{W}]=\frac{4}{n^2} \ee[D^{r-2}|W] \] for all $r$.

\bigskip

Lemma \ref{momcom2} performs some moment computations related to the pair
$(W,W')$.

\begin{lemma} \label{momcom2} Letting $D:=W'-W$, one has that:
\begin{enumerate}
\item $\ee(D|W)=-\frac{4}{n}(W-1)$.
\item $\ee(W)=1$.
\item $\ee(D^2|W)= \frac{8}{n}W - \frac{16}{n^2}(W-1)$.
\item $Var(W)=1$.
\item $\ee[D^4|W] = \frac{32}{n^2} \left( 2W^2 + \frac{12W - 8W^2}{n}+
\frac{8(1-W)}{n^2} \right)$.
\item $\ee[D^4|W] \leq \frac{256}{n^2}W^2 + \frac{256}{n^4} \ii
\{W=0\}$.
\end{enumerate}
\end{lemma}

\begin{proof} For part 1, by the construction of $(W,W')$ one has that
\begin{eqnarray*} \ee(D|i) & = & K(i,i+1) [W(i+1)-W(i)] + K(i,i-1) [W(i-1)-W(i)]\\
& = & \frac{(n+1-i)(n-2i) }{n(n+1-2i)} \left( \frac{4i}{n}-2 \right) +
\frac{i(n-2i+2)}{n(n+1-2i)} \left(2-\frac{4(i-1)}{n} \right).
\end{eqnarray*} Elementary simplifications show that this to equal
$-\frac{4}{n} (W(i)-1)$.

For part 2, since $W$ and $W'$ have the same law, one has that $\ee(D)=0$.
By part 1, \[ \ee(D) = \ee[\ee(D|W)] = - \frac{4}{n} \ee(W-1),
\] and the result follows.

For part 3, the construction of $(W,W')$ gives that
\begin{eqnarray*} \ee[D^2|i] & = & K(i,i+1) [W(i+1)-W(i)]^2 +
K(i,i-1) [W(i-1)-W(i)]^2 \\
& = & \frac{(n+1-i)(n-2i)\left( \frac{4i}{n}-2 \right)^2}{n(n+1-2i)} +
\frac{i(n-2i+2) \left(2-\frac{4(i-1)}{n} \right)^2}{n(n+1-2i)}.
\end{eqnarray*} Part 3 now follows by elementary algebra.

For part 4, observe that
\begin{eqnarray*}
\ee[D^2] & = & \ee[\ee[(W'-W)^2|W]]\\
& = & \ee[(W')^2] + \ee(W^2) - \ee[2W \ee(W'|W)]\\
& = & 2 \ee(W^2) - \ee[2W \ee(W'|W)]\\
& = & 2 \ee(W^2) - \ee \left[ 2W \left( (1-\frac{4}{n})W + \frac{4}{n} \right)
\right]\\
& = & \frac{8}{n} \ee(W^2) - \frac{8}{n}. \end{eqnarray*} The third
equality used that $W$ and $W'$ have the same distribution. The fourth
equality used part 1, and the final equality used part 2. Now parts 2 and
3 imply that $\ee[D^2]= \frac{8}{n}$. Thus $\ee(W^2)=2$, which together
with part 2 implies that $Var(W)=1$.

For part 5, note by  the construction of $(W,W')$ that
\begin{eqnarray*} \ee(D^4|i)
 & = & K(i,i+1) [W(i+1)-W(i)]^4 + K(i,i-1) [W(i-1)-W(i)]^4 \\
& = & \frac{(n+1-i)(n-2i)\left( \frac{4i}{n}-2 \right)^4}{n(n+1-2i)} +
\frac{i(n-2i+2) \left( 2 - \frac{4(i-1)}{n} \right)^4}{n(n+1-2i)}.
\end{eqnarray*} Elementary simplifications complete the proof of part 5.

Part 6 will follow from part 5. If $W=0$, then $\ee[D^4|W] =
\frac{256}{n^4}$, so part 6 is valid in this case. If $W \neq 0$, then by
the definition of $W$ it follows that $W \geq \frac{4}{n}$. Note that \[
\ee[D^4|W] - \frac{256}{n^2} W^2 = \frac{64}{n^2} \left( -3W^2 +
\frac{6W}{n} + \frac{4}{n^2} \right) -\frac{256 W^2}{n^3} - \frac{256
W}{n^4}.\] It is easy to see that $-3W^2+\frac{6W}{n}+\frac{4}{n^2}<0$ if
$W \geq \frac{4}{n}$, implying that $\ee[D^4|W] \leq \frac{256}{n^2} W^2$
if $W \neq 0$. \end{proof}

{\sc Proof of Proposition \ref{large}} One applies Theorem \ref{basicexp}
to the pair $(W,W')$. By Part 1 of Lemma \ref{momcom2}, the hypotheses are
satisfied with $\lambda=\frac{4}{n}$.

Consider the first error term in Theorem \ref{basicexp}. By parts 3 and 4
of Lemma \ref{momcom2}, \begin{eqnarray*}
\frac{\ee|2 \lambda W - \ee[D^2|W]|}{2 \lambda t} & = & \frac{2}{tn} \ee|W-1|\\
& \leq & \frac{2}{tn} \sqrt{\ee(W-1)^2}\\ & = & \frac{2}{tn}.
\end{eqnarray*}

Consider the second error term in Theorem \ref{basicexp}. By the
Cauchy-Schwarz inequality, \[ \ee|D|^3 \leq \sqrt{\ee[D^2] \ee[D^4]}.\]
Taking expectations in part 3 Lemma \ref{momcom2} gives that $\ee[D^2] =
\frac{8}{n}$. Taking expectations in part 5 of Lemma \ref{momcom2} gives
that $\ee[D^4] = \frac{128}{n^2}-\frac{128}{n^3} \leq \frac{128}{n^2}$.
Thus the second error term in Theorem \ref{basicexp} is at most $\frac{2
\max \{t^{-1},2t^{-2}\} }{\sqrt{n}}$.

To bound the third error term in Theorem \ref{basicexp}, one applies
Theorem \ref{tool3} with $\kappa=2$. Note from part 4 of Lemma
\ref{momcom2} that $\ee|W-1| \leq \sqrt{\ee(W-1)^2} = 1$. It is necessary
to bound \begin{align*} \epsilon_1(t) = \ee\bigl[\ee(D^2|W)\ii\{\ee(D^2|W)
> 2\lambda(W + t) \text{
 \textup{or} } \ee(D^4|W) > 4\lambda^2(\kappa^2W^2+\kappa^2t^2)\}\bigr].
\end{align*} Part 3 of Lemma \ref{momcom2} implies that $\ee[D^2|W)] > 2
\lambda(W+t)$ if and only if $(W-1)< -\frac{tn}{2}$. Part 5 of Lemma
\ref{momcom2} implies that
\[ \ee[D^4|W)] > 4 \lambda^2 (\kappa^2 W^2 + \kappa^2 t^2) \] can happen
only if $W=0$. Thus
\begin{equation} \label{ep1} \epsilon_1(t)  \leq \ee \left[ \ee[D^2|W]
\ii \left \{ W-1 < - \frac{tn}{2} \right\} \right]\\  + \ \pp(W=0)
\ee[D^2|W=0]. \end{equation} To bound the first term in (\ref{ep1}), note
by part 3 of Lemma \ref{momcom2} that
\[ \ee[D^2|W] = 2 \lambda + (2 \lambda - \lambda^2)(W-1).\] Since $n
\geq 2$, one has that $2 \lambda - \lambda^2 \geq 0$. It follows that if
$W-1<0$, then $\ee[D^2|W] \leq 2 \lambda$. Hence the first term in
(\ref{ep1}) is at most $\frac{8}{n} \pp( W-1< -\frac{tn}{2})$. By
Chebyshev's inequality, this is at most $\frac{32}{n^3 t^2}$.  To bound
the second term in (\ref{ep1}), note that $\pp(W=0) \leq \frac{2}{n}$.
Also part 3 of Lemma \ref{momcom2} gives that $\ee[D^2|W=0] =
\frac{16}{n^2}$, so that the second term in (\ref{ep1}) is at most
$\frac{32}{n^3}$. Summarizing, we have shown that $\epsilon_1(t) \leq
\frac{32}{n^3}(1+\frac{1}{t^2})$.

It is also necessary to bound \[ \epsilon_2(t) = \pp\{\ee(D^2|W) <
2\lambda(W - \frc{1}{4}t) \text{
  \textup{or}
} \ee(D^4|W)  > 4\lambda^2(\kappa^2W^2+\kappa^2t^2)\}. \] Part 3 of Lemma
\ref{momcom2} gives that \[ \ee[D^2|W] < 2 \lambda W - \frac{\lambda t}{2}
\] if and only if $(W-1) > \frac{tn}{8}$. Since $W$ has mean and
variance 1, Chebyshev's inequality implies that this occurs with
probability at most $\frac{64}{t^2n^2}$. By part 5 of Lemma \ref{momcom2},
\[ \ee[D^4|W)] > 4 \lambda^2 (\kappa^2 W^2 + \kappa^2 t^2) \] implies
that $W=0$. Since $\pp(W=0) \leq \frac{2}{n}$, it follows that
$\epsilon_2(t) \leq \frac{2}{n} + \frac{64}{t^2n^2}$.

Summarizing, the bounds on $\ee|W-1|, \epsilon_1(t),\epsilon_2(t)$ give
that the third error term in Theorem \ref{basicexp} is at most $\frac{B
\cdot \max \{1,t^{-3} \}}{\sqrt{n}}$ where $B$ is a universal constant.
Adding this to the first two error terms completes the proof.  \hfill
$\Box$

\subsection{Lower bound} \label{lower}

The purpose of this subsection is to prove the lower bound from Theorem
\ref{sharper} in the introduction.

\begin{prop} \label{lb} There is a sequence of $n$'s
tending to infinity, and corresponding $t_n$'s such that \[ |\pp(W \leq
t_n) - \pp(Z \leq t_n)| = \frac{2e^{-2}}{\sqrt{n}} + O(1/n).\]
\end{prop}

\begin{proof} Given $n$, define $i= \lceil \frac{n}{2} - \sqrt{n} \rceil$ and
$t_n=\frac{(n-2i)(n-2i+2)}{2n}$. The sequence of $n$'s will consist of
even perfect squares; then  $i= \lceil \frac{n}{2} - \sqrt{n} \rceil =
\frac{n}{2} - \sqrt{n}$ is an integer and the ceiling function can be
ignored.

Clearly \[ \pp(Z \geq t_n)=e^{(-2-\frac{2}{\sqrt{n}})}=e^{-2} \left(1 -
\frac{2}{\sqrt{n}}+O(\frac{1}{n}) \right).\] Also \[ \pp(W \geq t_n) =
\pp(i \leq \frac{n}{2} - \sqrt{n}) = \frac{{n \choose
\frac{n}{2}-\sqrt{n}}}{{n \choose \frac{n}{2}}}.\] Note that for integral
$a$,
\begin{eqnarray*}
\frac{{n \choose \frac{n}{2}-a}}{{n \choose \frac{n}{2}}} & = &
\frac{1}{(1+\frac{2a}{n})} \prod_{j=1}^{a-1}
\frac{(1-\frac{2j}{n})}{(1+\frac{2j}{n})}\\
& = & \frac{1}{(1+\frac{2a}{n})} e^{\sum_{j=1}^{a-1}
[\log(1-\frac{2j}{n}) - \log(1+\frac{2j}{n})]} \\
& = & \frac{1}{(1+\frac{2a}{n})} e^{\sum_{j=1}^{a-1} -\frac{4j}{n} +
O(\frac{j}{n})^3}. \end{eqnarray*} Since $a=\sqrt{n}$, one obtains that
\begin{eqnarray*}
\pp(W \geq t_n) & = & \frac{1}{(1+\frac{2}{\sqrt{n}})} e^{-2 +
\frac{2}{\sqrt{n}}+O(\frac{1}{n})}\\
& = & e^{-2} + O(1/n), \end{eqnarray*} and the result follows.
\end{proof}

{\it Remark:} Similar ideas give another proof of an $O(n^{-1/2})$ upper
bound for $|\pp(W \leq t)-\pp(Z \leq t)|$, when $t$ is fixed. This
argument was sketched to us by a referee of a much earlier (2006) version
of this paper, and goes as follows. The first step is to consider
$t=\frac{2(j^2+j)}{n}$ where $j$ is integral. Then $\pp(W \geq t)=
\frac{{n \choose \frac{n}{2}-j}}{{n \choose \frac{n}{2}}}$. From page 1077
of \cite{O}, one has the asymptotics \begin{equation} \label{as} \frac{{n
\choose \frac{n}{2}-j}}{{n \choose \frac{n}{2}}} =
e^{-\frac{2j^2}{n}+O(\frac{|j|^3}{n^2})} \end{equation} for $j \leq n/4$.
Since $t$ is fixed, one has that $j=O(n^{1/2})$ and so \begin{equation}
\label{24} |\pp(W \geq t) - \pp(Z \geq t)| =
|e^{-\frac{2j^2}{n}+O(\frac{|j|^3}{n^2})} - e^{\frac{-2(j^2+j)}{n}}| =
O(n^{-1/2}).\end{equation} The second step is to give a discretization
argument allowing one to also use non-integral $j$. The point is that for
fixed $t$ and $n$ growing, one can find an integer $j$ such that
$\frac{2(j^2+j)}{n} \leq t \leq \frac{2[(j+1)^2+(j+1)]}{n}$. For
$j=O(n^{1/2})$, one easily checks that \begin{equation} \label{25} \left|
\pp \left( Z \geq \frac{2(j^2+j)}{n} \right) - \pp \left( Z \geq
\frac{2[(j+1)^2+(j+1)]}{n} \right) \right| = O(n^{-1/2}) \end{equation}
and (using (\ref{as})) that \begin{equation} \label{26} \left| \pp \left(W
\geq \frac{2(j^2+j)}{n} \right)- \pp \left(W \geq
\frac{2[(j+1)^2+(j+1)]}{n} \right) \right| = \frac{{n \choose n/2-j}}{{n
\choose n/2}} - \frac{{n \choose n/2-j-1}}{{n \choose n/2}} = O(n^{-1/2}).
\end{equation} The $O(n^{-1/2})$ upper bound for $|\pp(W \leq t)-
\pp(Z \leq t)|$ with arbitrary $t>0$ fixed follows from \eqref{24},
\eqref{25}, and \eqref{26}.

\appendix
\section{Exchangeable Pair and Moment Computations: Algebraic Approach} \label{Alg}

The purpose of this appendix is to explain an algebraic approach to the
construction of the exchangeable pair $(W,W')$ in Subsection \ref{larget}
and to the moment computations in Lemma \ref{momcom2}. Since the
exchangeable pair in Subsection \ref{smallt} is related to that of
Subsection \ref{larget} (see the discussion in Subsection \ref{larget}),
this appendix gives insight into that exchangeable pair too. Throughout we
give results for the Johnson graph $J(n,k)$, as this contains the
Bernoulli-Laplace Markov chain as a special case $k=\frac{n}{2}$.

Let $G$ be a finite group and $K$ a subgroup of $G$. One calls $(G,K)$ a
Gelfand pair if the induced representation $1_K^G$ is multiplicity free.
For background on this concept, see Chapter 3 of \cite{D}, Chapter 7 of
\cite{Mac}, or Chapters 19 and 20 of \cite{T1}.

Suppose that $(G,K)$ is a Gelfand pair, so that $1_K^G$ decomposes as
$\bigoplus_{i=0}^s V_i$, where $V_0$ is the trivial module. Letting $d_i$
be the dimension of $V_i$, one can define a probability measure $\pi$ on
$\{0,\cdots,s \}$ by $\pi(i) = \frac{d_i}{|G/K|}$. Associated to each
value of $i$ between $0$ and $s$ is a ``spherical function'' $\omega_i$,
which is a certain map from the double cosets of $K$ in $G$ to the complex
numbers. Hence $\pi$ can be viewed as a probability measure on spherical
functions.

    The spectrum of the Johnson graph $J(n,k)$ can be understood
    in the language of spherical functions of Gelfand pairs; this
    goes back to \cite{DS}, which used this viewpoint to study the
    convergence rate of random walk on $J(n,k)$. To describe this,
    suppose without loss of generality that $0 \leq k \leq
    \frac{n}{2}$. Let $G$ be the symmetric group $S_n$, and $K$
    the subgroup $S_k \times S_{n-k}$. Then the space $G/K$ is in
    bijection with the vertices of the Johnson graph. There are
    $k+1$ spherical functions $\{\omega_0,\cdots,\omega_{k} \}$,
    and the dimension $d_i$ is equal to ${n \choose i} - {n
    \choose i-1}$ if $1 \leq i \leq k$ and to 1 if $i=0$. The
    double cosets $K_0,K_1,\cdots, K_k$ of $K$ in $G$ are also
    indexed by the numbers $0,1,\cdots,k$; the double coset
    corresponding to $j$ consists of those permutations $\tau$ in
    $S_n$ such that $|\{1,\cdots,k\} \cap
    \{\tau(1),\cdots,\tau(k)\}| = k-j$. Letting $\omega_i(j)$
    denote the value of $\omega_i$ on the double coset indexed by
    $j$, it is known that \[ \omega_i(j) = \sum_{m=0}^i
    \frac{(-i)_m (i-n-1)_m (-j)_m}{(k-n)_m (-k)_m m!} \] where
    $(j)_m = j(j+1) \cdots (j+m-1)$ for $m \geq 1$ and
    $(j)_0=1$. The spectrum of random walk on the Johnson graph
    consists of the numbers $\omega_i(1)$ with multiplicity $d_i$.

Specializing to $k=\frac{n}{2}$ in the previous paragraph, the random
variable $W$ studied in Section \ref{n/2} is equal to $W(i) = \frac{n}{2}
\omega_i(1) + 1$, so up to constants is a random spherical function of the
Gelfand pair $(G,K)$. Section 4 of the paper \cite{F4} used Stein's method
to study random spherical functions of Gelfand pairs. Although the
examples studied there were all for normal approximation, many of the
results are general. For example, an exchangeable pair $(W,W')$ was
constructed using a reversible Markov chain. Specializing to the Gelfand
pair corresponding to $J(n,k)$, the Markov chain is on the set
$\{0,1,\cdots,k\}$ and transitions from $i$ to $j$ with probability \[
L(i,j) := \frac{d_j}{|G|} \sum_{r=0}^k |K_r| \omega_i(K_r) \omega_1(K_r)
\overline{\omega_j(K_r)}.\]

Proposition \ref{orthog} proves that the Markov chain $L$ is a birth-death
chain (and specializes to the birth-death chain of Subsection \ref{larget}
when $k=\frac{n}{2}$). This is interesting, since from the definition of
$L$ it is not even evident that it is a birth-death chain.

\begin{prop} \label{orthog} The Markov chain $L$ on the set $\{0,1,\cdots,k \}$
is a birth-death chain with transition probabilities
\[ L(i,i+1) = \frac{n(n+1-i)(n-i-k)(k-i)}{k(n-k)(n+1-2i)(n-2i)} \]
\[L(i,i-1) = \frac{in(n+1-i-k)(k+1-i)}{k(n-k)(n+2-2i)(n+1-2i)} \]
 \[ L(i,i) = \frac{i(n+1-i)(n-2k)^2}{k(n-k)(n-2i)(n+2-2i)}\]
\end{prop}

\begin{proof} The spherical function $\omega_i(K_r)$ is the Hahn polynomial
$Q_n(x;\alpha,\beta,N)$ where $x=r,n=i,N=k,\alpha=k-n-1,\beta=-k-1$.
Properties of these polynomials are given on pages 33-34 of \cite{KoSw}.
In particular, they satisfy a recurrence relation \[ -r \omega_i(K_r) =
A_i \omega_{i+1}(K_r) - (A_i+B_i) \omega_i(K_r) + B_i \omega_{i-1}(K_r) \]
where \[ A_i = \frac{(n+1-i)(n-k-i)(k-i)}{(n+1-2i)(n-2i)} \] and \[ B_i =
\frac{i(n+1-k-i)(k+1-i)}{(n+2-2i)(n+1-2i)}.\]

Since $\omega_1(K_r) = 1 - \frac{nr}{k(n-k)}$, it follows that
\begin{eqnarray*} \omega_1(K_r) \omega_i(K_r) & = &
 \frac{n(n+1-i)(n-i-k)(k-i)}{k(n-k)(n+1-2i)(n-2i)} \omega_{i+1}(K_r)\\
& & + \frac{i(n+1-i)(n-2k)^2}{k(n-k)(n-2i)(n+2-2i)} \omega_i(K_r)\\
& & + \frac{in(n+1-i-k)(k+1-i)}{k(n-k)(n+2-2i)(n+1-2i)} \omega_{i-1}(K_r).
\end{eqnarray*} The result now follows immediately from the orthogonality
relations for Hahn-polynomials \cite{KoSw}, which are a special case of
the orthogonality relations for spherical functions of a Gelfand pair
\cite {Mac}. \end{proof}

    To conclude, we note that there is an algebraic way to compute
    the moments $\ee(W'-W)^m$ and the conditional moments
    $\ee[(W'-W)^m|i]$. The interesting point about this approach
    is that it does not require one to explicitly compute the
    transition probabilities of the Markov chain $L$, or even to
    know that in this particular case it is a birth-death
    chain. Moreover, some of the quantities which appear have direct
    interpretations in terms of random walk on the Johnson graph.

    To be precise, Lemma 4.12 of \cite{F4} implies that
    $\ee(W'-W)^m$ is equal to \[ \left( \frac{|K_1|}{|K|}
    \right)^{m/2} \sum_{l=0}^m (-1)^{m-l} {m \choose l}
    \sum_{r=0}^s \frac{|K|}{|K_r|} \omega_1(K_r) p_l(K_r)
    p_{m-l}(K_r). \] Here $p_j(K_r)$ is the chance that random
    walk on the Johnson graph $J(n,k)$ started at a particular
    vertex, is distance r away from the start vertex after j
    steps. Also, the proof of the lemma gives that
    $\ee[(W'-W)^m|i]$ is equal to \[ \left( \frac{|K_1|}{|K|}
    \right)^{m/2} \sum_{l=0}^m (-1)^{m-l} {m \choose l}
    \omega_i(K_1)^{m-l} \sum_{r=0}^s \omega_i(K_r) \omega_1(K_r)
    p_l(K_r). \] These expressions are easily evaluated for small
    $m$, and one obtains another proof of Lemma \ref{momcom2}.

\section*{Acknowledgements} Fulman was supported by NSA grant H98230-05-1-0031
and NSF grants DMS-0503901 and DMS-0802082.


\begin{thebibliography}{AAA}

\bibitem [CS]{CS} Chen, L. and Shao, Q., Stein's method for normal
approximation, in {\it An introduction to Stein's method}, Lecture
Notes Series, Institute for Mathematical Sciences, National University
of Singapore, Volume 4, 2005, 1-59.

\bibitem [D]{D} Diaconis, P., {\it Group representations in
probability and statistics}, Institute of Mathematical Statistics
Lecture Notes, Volume 11, 1988.

\bibitem [DS]{DS} Diaconis, P. and Shahshahani, M., Time to reach
stationarity in the Bernoulli-Laplace diffusion model, {\it Siam J. of
Math. Anal.} {\bf 18} (1987), 208-218.

\bibitem [F1]{F1} Fulman, J., Stein's method and Plancherel measure of
the symmetric group, {\it Transac. Amer. Math. Soc.} {\bf 357} (2005),
555-570.

\bibitem [F2]{F2} Fulman, J., Martingales and character ratios,
{\it Trans. Amer. Math. Soc.} {\bf 358} (2006), 4533-4552.

\bibitem [F3]{F3} Fulman, J., An inductive proof of the Berry-Ess\'{e}en theorem
for character ratios, {\it Ann. Comb.} {\bf 10} (2006), 319-332.

\bibitem [F4]{F4} Fulman, J., Stein's method and random character
ratios, {\it Trans. Amer. Math. Soc.} {\bf 360} (2008), 3687-3730.

\bibitem [GT]{GT} G\"{o}tze, F. and Tikhomirov, A. N.,
Limit theorems for spectra of random matrices with martingale structure,
{\it Theory Probab. Appl.} {\bf 51} (2007), 42-64.

\bibitem [Ho1]{Ho1} Hora, A., Central limit theorems and asymptotic
spectral analysis on large graphs, {\it
Inf. Dim. Anal. Quant. Prob. and Rel. Topics} {\bf 1} (1998), 221-246.

\bibitem [Ho2]{Ho2} Hora, A., Central limit theorem for the adjacency operators on the
infinite symmetric group, {\it Comm. Math. Phys.} {\bf 195} (1998),
405-416.

\bibitem [IO]{IO} Ivanov, V. and Olshanski, G., Kerov's central limit theorem
for the Plancherel measure on Young diagrams, in {\it Symmetric functions
2001: surveys of developments and perspectives}, 93-151, NATO Sci. Ser. II
Math. Phys. Chem., 74, Kluwer Acad. Publ., Dordrecht, 2002.

\bibitem [Ke]{Ke} Kerov, S., Gaussian limit for the Plancherel measure of
the symmetric group, {\it Compt. Rend. Acad. Sci. Paris, Serie I}, {\bf
316} (1993), 303-308.

\bibitem [KoSw]{KoSw} Koekoek, R. and Swarttouw, R., The Askey-scheme
of hypergeometric orthogonal polynomials and its q-analog,
arXiv:math.CA/9602214 (1996).

\bibitem [Lu]{Lu} Luk, H.M., Stein's method for the gamma distribution and related
statistical applications, Ph.D. thesis, University of Southern California,
1994.

\bibitem [Mc]{Mac} Macdonald, I., {\it Symmetric functions and Hall polynomials},
Second edition, Oxford University Press, New York, 1995.

\bibitem [Mn]{Man} Mann, B., Stein's method for $\chi^2$ of a
multinomial, unpublished manuscript (1997).

\bibitem [O]{O} Odlyzko, A., Asymptotic enumeration methods, in
{\it Handbook of combinatorics, Vol. 2}, 1063-1229, Elsevier, Amsterdam,
1995.

\bibitem [Re]{Re} Reinert, G., Three general approaches to Stein's method,
in {\it An introduction to Stein's method}, Lecture Notes Series,
Institute for Mathematical Sciences, National University of Singapore,
Volume 4 (1994), 183-221.

\bibitem [RR]{RR} Rinott, Y. and Rotar, V., On coupling constructions and rates
in the CLT for dependent summands with applications to the antivoter model
and weighted $U$-statistics. {\it Ann. Appl. Probab.} {\bf 7} (1997),
1080-1105.

\bibitem [Ro]{Ro} R\"{o}llin, A., A note on the exchangeability condition in Stein's
method, {\it Statist. Probab. Lett.}, in press (2008).

\bibitem [Sc]{Sc} Scarabotti, F., Time to reach stationarity in the
Bernoulli-Laplace diffusion model with many urns, {\it Adv. Appl. Math.}
{\bf 18} (1997), 351-371.

\bibitem [ShSu]{ShSu} Shao, Q., and Su, Z., The Berry-Ess\'{e}en bound for
character ratios, {\it Proc. Amer. Math. Soc.} {\bf 134} (2006),
2153-2159.

\bibitem[Sn]{Sn} Sniady, P., Gaussian fluctuations of characters of symmetric
groups and of Young diagrams, {\it Probab. Theory Related Fields} {\bf
136} (2006), 263-297.

\bibitem[St1]{St1} Stein, C., {\it Approximate computation of expectations},
Institute of Mathematical Statistics Lecture Notes, Volume 7, 1986.

\bibitem[St2]{St2} Stein, C., with Diaconis, P., Holmes, S., and Reinert,
G., Use of exchangeable pairs in the analysis of simulations, in {\it
Stein's method: expository lectures and applications},  69-77, IMS Lecture
Notes Monogr. Ser., Volume 46, 2004.

\bibitem[T1]{T1} Terras, A., {\it Fourier analysis on finite groups and
applications}, London Mathematical Society Student Texts 43, Cambridge
University Press, Cambridge, 1999.

\bibitem[T2]{T2} Terras, A., Survey of spectra of Laplacians on finite symmetric spaces,
{\it Experimental Math.}, {\bf 5} (1996), 15-32.

\end{thebibliography}
\end{document}